\documentclass[12pt]{article}
\usepackage{mathrsfs}
\usepackage{amssymb}
\usepackage{amsfonts}
\usepackage{latexsym}
\usepackage{amsthm}
\usepackage{pictex}

\usepackage[T2A]{fontenc}     % внутренняя {{T2A}} кодировка {{TeX}}
\usepackage[cp1251]{inputenc} % кодировка - можно использовать {{[cp866] [koi8-r]}}
\usepackage{graphicx}

\newcommand{\er}[1]{{\rm(\ref{#1})}}
\def\lb{\label}
\theoremstyle{plain}
\newtheorem{theorem}{\bf Theorem}[section]
\newtheorem{lemma}[theorem]{\bf Lemma}
\theoremstyle{remark}

\newtheorem{proposition}[theorem]{\bf Proposition}

\begin{document}

\def\a{\alpha} \def\cA{{\cal A}} \def\bA{{\bf A}}  \def\mA{{\mathscr A}}
\def\b{\beta}  \def\cB{{\cal B}} \def\bB{{\bf B}}  \def\mB{{\mathscr B}}
\def\g{\gamma} \def\cC{{\cal C}} \def\bC{{\bf C}}  \def\mC{{\mathscr C}}
\def\G{\Gamma} \def\cD{{\cal D}} \def\bD{{\bf D}}  \def\mD{{\mathscr D}}
\def\d{\delta} \def\cE{{\cal E}} \def\bE{{\bf E}}  \def\mE{{\mathscr E}}
\def\D{\Delta} \def\cF{{\cal F}} \def\bF{{\bf F}}  \def\mF{{\mathscr F}}
\def\c{\chi}   \def\cG{{\cal G}} \def\bG{{\bf G}}  \def\mG{{\mathscr G}}
\def\z{\zeta}  \def\cH{{\cal H}} \def\bH{{\bf H}}  \def\mH{{\mathscr H}}
\def\e{\eta}   \def\cI{{\cal I}} \def\bI{{\bf I}}  \def\mI{{\mathscr I}}
\def\p{\psi}   \def\cJ{{\cal J}} \def\bJ{{\bf J}}  \def\mJ{{\mathscr J}}
\def\vT{\Theta}\def\cK{{\cal K}} \def\bK{{\bf K}}  \def\mK{{\mathscr K}}
\def\k{\kappa} \def\cL{{\cal L}} \def\bL{{\bf L}}  \def\mL{{\mathscr L}}
\def\l{\lambda}\def\cM{{\cal M}} \def\bM{{\bf M}}  \def\mM{{\mathscr M}}
\def\L{\Lambda}\def\cN{{\cal N}} \def\bN{{\bf N}}  \def\mN{{\mathscr N}}
\def\m{\mu}    \def\cO{{\cal O}} \def\bO{{\bf O}}  \def\mO{{\mathscr O}}
\def\n{\nu}    \def\cP{{\cal P}} \def\bP{{\bf P}}  \def\mP{{\mathscr P}}
\def\r{\rho}   \def\cQ{{\cal Q}} \def\bQ{{\bf Q}}  \def\mQ{{\mathscr Q}}
\def\s{\sigma} \def\cR{{\cal R}} \def\bR{{\bf R}}  \def\mR{{\mathscr R}}
\def\S{\Sigma} \def\cS{{\cal S}} \def\bS{{\bf S}}  \def\mS{{\mathscr S}}
\def\t{\tau}   \def\cT{{\cal T}} \def\bT{{\bf T}}  \def\mT{{\mathscr T}}
\def\f{\phi}   \def\cU{{\cal U}} \def\bU{{\bf U}}  \def\mU{{\mathscr U}}
\def\F{\Phi}   \def\cV{{\cal V}} \def\bV{{\bf V}}  \def\mV{{\mathscr V}}
\def\P{\Psi}   \def\cW{{\cal W}} \def\bW{{\bf W}}  \def\mW{{\mathscr W}}
\def\o{\omega} \def\cX{{\cal X}} \def\bX{{\bf X}}  \def\mX{{\mathscr X}}
\def\x{\xi}    \def\cY{{\cal Y}} \def\bY{{\bf Y}}  \def\mY{{\mathscr Y}}
\def\X{\Xi}    \def\cZ{{\cal Z}} \def\bZ{{\bf Z}}  \def\mZ{{\mathscr Z}}
\def\O{\Omega}
\def\ve{\varepsilon}
\def\vt{\vartheta}
\def\vp{\varphi}
\def\vk{\varkappa}

%%%%%%%%%%%%%%%%%%%%%
\def\mM{M}
\def\mB{B}
\def\mR{R}

\def\mA{{\mathscr A}}
\def\mB{{\mathscr B}}
\def\mC{{\mathscr C}}
\def\mD{{\mathscr D}}
\def\mE{{\mathscr E}}
\def\mF{{\mathscr F}}
\def\mG{{\mathscr G}}
\def\mH{{\mathscr H}}
\def\mI{{\mathscr I}}
\def\mJ{{\mathscr J}}
\def\mK{{\mathscr K}}
\def\mL{{\mathscr L}}
\def\mM{{\mathscr M}}
\def\mN{{\mathscr N}}
\def\mO{{\mathscr O}}
\def\mP{{\mathscr P}}
\def\mQ{{\mathscr Q}}
\def\mR{{\mathscr R}}
\def\mS{{\mathscr S}}
\def\mT{{\mathscr T}}
\def\mU{{\mathscr U}}
\def\mV{{\mathscr V}}
\def\mW{{\mathscr W}}
\def\mX{{\mathscr X}}
\def\mY{{\mathscr Y}}
\def\mZ{{\mathscr Z}}

\def\Z{{\Bbb Z}}
\def\R{{\Bbb R}}
\def\C{{\Bbb C}}
\def\T{{\Bbb T}}
\def\N{{\Bbb N}}
\def\S{{\Bbb S}}
\def\H{{\Bbb H}}

\def\qqq{\qquad}
\def\qq{\quad}
\newcommand{\ma}{\begin{pmatrix}}
\newcommand{\am}{\end{pmatrix}}
\let\ge\geqslant
\let\le\leqslant
\let\geq\geqslant
\let\leq\leqslant
\def\ma{\left(\begin{array}{cc}}
\def\am{\end{array}\right)}
\def\iint{\int\!\!\!\int}
\def\lt{\biggl}
\def\rt{\biggr}
\let\geq\geqslant
\let\leq\leqslant
\def\[{\begin{equation}}
\def\]{\end{equation}}
\def\wt{\widetilde}
\def\pa{\partial}
\def\sm{\setminus}
\def\es{\emptyset}
\def\no{\noindent}
\def\ol{\overline}
\def\iy{\infty}
\def\ev{\equiv}
\def\/{\over}
\def\ts{\times}
\def\os{\oplus}
\def\ss{\subset}
\def\h{\hat}
\def\Re{\mathop{\rm Re}\nolimits}
\def\Im{\mathop{\rm Im}\nolimits}
\def\supp{\mathop{\rm supp}\nolimits}
\def\sign{\mathop{\rm sign}\nolimits}
\def\Ran{\mathop{\rm Ran}\nolimits}
\def\Ker{\mathop{\rm Ker}\nolimits}
\def\Tr{\mathop{\rm Tr}\nolimits}
\def\const{\mathop{\rm const}\nolimits}
\def\dist{\mathop{\rm dist}\nolimits}
\def\diag{\mathop{\rm diag}\nolimits}
\def\Wr{\mathop{\rm Wr}\nolimits}
\def\BBox{\hspace{1mm}\vrule height6pt width5.5pt depth0pt \hspace{6pt}}

\def\Diag{\mathop{\rm Diag}\nolimits}

\def\Twelve{
\font\Tenmsa=msam10 scaled 1200 \font\Sevenmsa=msam7 scaled 1200
\font\Fivemsa=msam5 scaled 1200 \textfont\msbfam=\Tenmsb
\scriptfont\msbfam=\Sevenmsb \scriptscriptfont\msbfam=\Fivemsb

\font\Teneufm=eufm10 scaled 1200 \font\Seveneufm=eufm7 scaled 1200
\font\Fiveeufm=eufm5 scaled 1200
%\newfam\eufmfam
\textfont\eufmfam=\Teneufm \scriptfont\eufmfam=\Seveneufm
\scriptscriptfont\eufmfam=\Fiveeufm}

\def\Ten{
\textfont\msafam=\tenmsa \scriptfont\msafam=\sevenmsa
\scriptscriptfont\msafam=\fivemsa

\textfont\msbfam=\tenmsb \scriptfont\msbfam=\sevenmsb
\scriptscriptfont\msbfam=\fivemsb

\textfont\eufmfam=\teneufm \scriptfont\eufmfam=\seveneufm
\scriptscriptfont\eufmfam=\fiveeufm}

\title {Lyapunov functions for periodic matrix-valued Jacobi operators}

\author{Evgeny Korotyaev
\begin{footnote}
{ Institut f\"ur  Mathematik,  Humboldt Universit\"at zu Berlin,
Rudower Chaussee 25, 12489, Berlin, Germany, e-mail:
evgeny@math.hu-berlin.de}
\end{footnote}
 \and Anton Kutsenko
\begin{footnote}
{ Department of
 Mathematics of Sankt-Petersburg State University, Russia e-mail: kucenkoa@rambler.u}
\end{footnote}
}

\maketitle

\begin{abstract}
\no We consider  periodic matrix-valued Jacobi
operators. The spectrum of this operator is absolutely continuous
and consists of intervals separated by gaps. We define the
Lyapunov function, which is analytic on an associated
Riemann surface. On each sheet the Lyapunov function has the
standard properties of the Lyapunov function for the scalar case.
We show that this function has (real or complex) branch points, which we call resonances. We prove that there exist two types of gaps: i) stable
gaps, i.e., the endpoints are periodic and anti-periodic
eigenvalues, ii) unstable (resonance) gaps, i.e., the endpoints
are resonances (real branch points).
We show that some spectral data determine the spectrum (counting multiplicity) of the Jacobi operator.

\end{abstract}

%  \vskip 0.25cm

\section {Introduction and main results}
\setcounter{equation}{0}

\no Consider self-adjoint matrix-valued Jacobi operators
$\cJ$ acting on $\ell^2(\Z)^m$ and are given by
\[
 \lb{000}
 (\cJ y)_n=a_n y_{n+1}+b_ny_n+a_{n-1}^\top y_{n-1},\qq n\in\Z,\qq
 y_n\in \C^m, y=(y_n)_{n\in \Z}\in \ell^2(\Z)^m,m\ge 1,
\]
 where $a_n,b_n=b_n^\top$ are p-periodic sequences of the $m\ts m$ real matrices and $ \det a_n\ne 0$ for all $n\in\Z$. It is well known that the spectrum $\s(\cJ)$ of $\cJ$ is absolutely
continuous  and consists of non-degenerated intervals $[\l^+_{n-1},\l^-_n], \l^+_{n-1}<\l^-_n \le  \l^+_{n}, n=1,...,N<\iy$. %and let $\l^-_{N}=\l_0^-$.
These intervals are separated by the gaps $\g_n=(\l_n^-,\l_n^+)$ with the length $>0$. Introduce the fundamental $m\ts m$ matrix-valued solutions
$\vp=(\vp_n(z))_{n\in\Z}, \vt=(\vt_n(z))_{n\in\Z}$ of the
equation
\[
 \lb{001}
 a_ny_{n+1}+b_ny_n+a_{n-1}^\top y_{n-1}=z y_n,\qqq
\vp_{0}\ev \vt_1\ev 0,\ \vp_1\ev \vt_{0}\ev I_m,\qq
 (z,n)\in\C\ts\Z,
\]
where $I_m$ is the identity $m\ts m$ matrix. We
define the monodromy $2m\ts2m$ matrix $\mM_p$ and the trace
$T_n, n\ge 1$ by
\[
 \lb{002}
 \mM_p(z)=\ma \vt_{p}(z) & \vp_{p}(z) \\
           \vt_{p+1}(z) & \vp_{p+1}(z)
       \am,
       \qq T_n(z)=\Tr \mM_p^n(z).
\]
 Introduce the modified monodromy matrix $M$
and the determinant $D$ by
\[
 \lb{002b}
 M=P_0\mM_p P_0^{-1},\ \ \ P_0=\ma a_0& 0 \\ 0 & I_m \am,
 \qqq D(z,\t)=\det (M(z)-\t I_{2m}),\  \t,z\in\C.
\]
Let $\t_j=\t_j(z), j\in \N_{2m}$ be the eigenvalues of $M(z)$, where
$\N_{m}=\{1,,..,m\}$. An
eigenvalue of $M(z)$ is called a {\it multiplier}.
%It is a zero
%of the algebraic equation $D(z,\t)=0,\t\in\C$ for fixed $z\in\C$.

{\no\bf Theorem (Lyapunov-Poincar\'e)}
{\it  i) The following identities hold
\[
 \lb{700}
 M^\top JM=J,\ \ \ J=\ma 0 & I_m \\ -I_m & 0 \am,
\]
\[
\lb{701}
D(z,\t)=\t^{2m}D(z,\t^{-1}), \qq \  all \qq z,\t\in\C,
\t\ne 0.
\]
\no ii) If  $\t(z)$ is a multiplier of
multiplicity $d\ge 1$ for some $z\in\C$ (or $z\in\R$), then $\t^{-1}(z)$ (or $\ol\t(z)$) is a multiplier of multiplicity $d$. Moreover, each $M(z), z\in\C$,
has exactly $2m$ multipliers $\t_j^{\pm 1}(z), j\in \N_m$
and $\s(\cJ)=\cup_{j=1}^m\{z\in \C: |\t_j(z)|=1 \}$.

\no iii) If $\t(z)$ is a simple multiplier for some $z\in \C$ and $|\t(z)|=1$,
then $\t'(z)\not=0$.
%\no {{iv)}} For any $z,\t\in\C$, $\t\not=0$, the following
%identity holds
%\[
% zb{011}
% D(z,\t)=\a^{-1}{\t^m}z^{mN}+O(z^{mN-1}),\ \ \l\to\infty,
% \ \ \ \a=\prod_{j=1}^m\det A_j.
%\]
}
%It is well known that $ D(\t,\l)=\sum_0^{2m}\x_j(\l)\t^{2m-j}$,
%where the functions $\x_j$ are given by
%\[
%\lb{15} \x_0=1,\ \ \ \x_1=-2mT_1,\ \  \x_2=-{2m\/2}(T_2+T_1\x_1),\
%\ .. ... , \x_j=-{2m\/j}\sum_0^{j-1}T_{j-i}\x_i,..
%\]
%see p.331-333 [RS]. Using the identity \er{TL-2} we obtain
%\[
%\lb{poD}
% D(\t,\cdot)=(\t^{2m}+1)+\x_1(\t^{2m-1}+\t)+
%...+\x_{m-1}(\t^{m+1}+\t^{m-1})+\x_{m}\t^m.
%\]

 It is well known that the spectrum of the scalar ($m=1$) Jacobi operator $\cJ$ is absolutely continuous and $\s(\cJ)=\cup_1^p[\l_{n-1}^+,\l_{n}^-]$, where
$\l_0^{+}<\l_1^-\le \l_1^+< ... < \l_p^- \le \l_{p-1}^+ < \l_{p}^{-}$
and $\l_n^{\pm}$ are 2-periodic eigenvalues.
The intervals $[\l_{n-1}^+,\l_{n}^-],n\in \N_p$ are separated by
gaps $\g_n=(\l_n^-,\l_n^-)$ of lengths $|\g_n|\ge 0$. If a
gap $\g_n$ is degenerate, i.e. $|\g_n|=0$, then the corresponding
segments $\s_n$, $\s_{n+1}$ merge.  Recall that the Lyapunov function
$\D(z)={\vp_{p+1}(z)+\vt_p(z)\/2}$ and the spectrum $\s(\cJ)=\{\l\in\C:\ \D(z)\in [-1,1]\}$. Note that
$(-1)^{p-n}\D(\l_n^{\pm})=1$ for all $n\in \N_p$.
We recall well-known facts (see e.g. \cite{P} or \cite{KKu}).

\begin{theorem}\lb{16001}
A real polynomial $F$ is  the Lyapunov function for some scalar
$p$-periodic Jacobi operator with numbers $a_n>0,b_n\in\R,n\in\Z$
iff $F(z)=c z^p+O(z^{p-1})$ as $z\to \iy$ for some $c>0$ and
$
 F'(z_j)=0,\ \ (-1)^{p-j}F(z_j)\ge 1$ for all $j\in \N_{p-1}$
 and
for some $z_1<...<z_{p-1}$.
\end{theorem}

\no {\bf Remark.} 1) Zeros of $\D-t$ for any fixed $t\in [-1,1]$ and a constant $c>0$ determine $\D$.

\no 2)   $\s(\cJ)=\s(\wt\cJ)$ for some
Jacobi operators $\cJ,\wt\cJ$ iff $\D=\wt\D$, where $\D, \wt\D$
are the corresponding Lyapunov functions.

\no 3) If $a_n=1,b_n=0$ for all $n\in\Z$, then $\D=\mT_p({z\/2})$,
where $\mT_p$ is the Chebyshev polynomial.

Before we describe the content of our paper, we briefly
comment on background literature for matrix-valued Jacobi operators.
Inverse spectral theory for scalar periodic Jacobi operators
is well understood, see a book \cite{T} and papers
\cite{BGGK},
%\cite{HK},
\cite{K}, \cite{KKu}, \cite{P}, \cite{vM}
and references therein.
The corresponding theory for periodic Jacobi matrices with matrix-valued coefficients is still largely a wide
open field. Some new particular results were recently obtained  by
Gesztesy  and coauthors \cite{CG},\cite{CGR},\cite{GKM}
and see references therein.
Note that for finite Jacobi matrices with matrix-valued coefficients
the complete  solution of the inverse problem was given recently  by Chelkak and Korotyaev \cite{CK}.

The eigenvalues of $M(z)$ are the zeros of the equation $D(\t,z)=0$. This
is an algebraic equation in $\t$ of degree $2m$, where the coefficients
are polynomials of $z$. It is well known (see e.g. Chapter 8,
[Fo]) that the roots $\t_j(z),j\in \N_{2m}$ constitute one or several
branches of one or several analytic functions that have only
algebraic singularities in $\C$. Thus the number of eigenvalues of
$M(z)$ is a constant $N_e$ with the exception of some special
values of $z$ (see below the definition of a resonance). There is a finite number of such points on the plane.
If the functions $\t_j(z),j\in\N_{2m}$ are all distinct, then
$N_e=2m$. If some of them are identical, then $N_e<2m$ and
$M(z)$ is permanently degenerate.

If $m=1$, then the Riemann surface for the multipliers
 has 2 sheets, but the Lyapunov function is entire. Similarly, in the case $m\ge 2$ it is more convenient for us to  construct the Riemann surface for the Lyapunov functions given by
$\D_j(z)={1\/2}(\t_j(z)+\t_j^{-1}(z)), \ j\in\N_m
$, which has $m$ sheets (see the equation \er{T1-1}).
Note that the Lyapunov-Poincar\'e Theorem  gives $\s(\cJ)=\bigcup_{j=1}^m\{z: \D_j(z)\in [-1,1]\}$.
Let $\# A$ be a number of elements of a set $A$.

\begin{theorem} \lb{T1}
There exist analytic functions $\wt\D_s, s=1,..,p_0\le m$ on the
$p_s$-sheeted Riemann surface $\mR_s, p_s\ge 1$ having the
following properties:

\no i) There exist disjoint subsets $\o_s\ss \N_m,
s\in \N_{p_0}, \bigcup \o_s=\N_m$ such that all branches of
$\wt\D_s,s\in\N_{p_0}$ have the form $\D_j(z)={1\/2}(\t_j(z)+\t_j^{-1}(z)), \ j\in \o_s$ and satisfy
\[
\lb{T1-1}
{D(z,\t)\/(2\t)^m}=\prod_1^{p_0} \F_s(z,\n),\qq
\ \F_s(z,\n)=\prod_{j\in \o_s}(\n-\D_j(z)),\qq
\n={\t+\t^{-1}\/2},\qq \t\ne 0,
\]
for any $z,\t\in \C$, where the functions $\F_s(z,\n)$ are some polynomials with
respect to $\n,z\in\C$. If $\D_i=\D_j$ for some $i\in \o_k,
j\in \o_s$, then $\F_k=\F_s$ and $\wt \D_k=\wt \D_s$.

\no ii) Let some branch $\D_j,j\in \N_m$ be real analytic
on some interval $Y=(\a,\b)\ss\R$ and $-1<\D_j(z)<1$ for any
$z\in Y$. Then $\D_j'(z)\ne 0$ for each $z\in Y$ .

\no iii) All functions $\r, \r_s,s\in\N_{p_0}, \#\o_s\ge 2$ given by
\er{T1-2} are polynomials and
\[
\lb{T1-2} \r=\prod_{1}^{N_0}\r_s,\ \ \
\r_s(z)=\!\!\!\!\prod_{i<j, i,j\in \o_s}\!\!\!\!
(\D_i(z)-\D_j(z))^2,\ \ z\in\C,
\qq where \qq \r_s=1, \qq if \ \#\o_s=1.
\]
\no iv) Each endpoint of a gap $(\l_n^-,\l_n^+)$ is a periodic (or anti-periodic) eigenvalue or a real branch point of some $\D_j, j\in\N_{m}$, which is a zero of $\r$ (below such points are called resonances).

\no v) Let $\t_j^0, j\in\N_m$ be eigenvalues of a matrix
$A_p=(a_1a_2..a_p)^{-1}$. Then following asymptotics hold
\[
\lb{T1-3}
 \D_j(z)={z^p\/2}(\t_j^0+O(z^{-{1\/m}})),\ \ \
 \t_j(z)=z^p\t_j^0+O(z^{-{1\/m}})),\qq j\in\N_m,\ \ \
\]
\[
\lb{T1-4}
 \r_s(z)=z^{p|\o_s|(|\o_s|-1)}(c_s+O(z^{-{1\/m}})),\ \
 c_s=2^{|\o_s|(1-|\o_s|)}\prod_{k<j, k,j\in \o_s}(\t_k^0-\t_j^0)^2,\ \ s\in \N_{p_0}
\]
as $|z|\to\iy$ Moreover, if $\t_j^0\not=\t_k^0$ for all $j,k\in \o_s,
j\ne k$, then $c_s\not=0$ and $\r_s(z)=z^{|\o_s|(|\o_s|-1)}(c_s+O(z^{-1}))$ as $|z|\to\iy$.
\end{theorem}

  \vskip 0.25cm
{\bf Definition.} {\it The number $z_0$ is a {\bf resonance} of $\cJ$, if $z_0$ is a zero of $\r$ given by \er{T1-2}.}

\no {\bf Example.} Let  $a_n=\diag\{a_{1,n},a_{2,n},...,a_{m,n},\}$
and $b_n=\diag\{b_{1,n},b_{2,n},...,b_{m,n},\}$ for all $n\in \Z$
and for some $a_{j,n} >0, b_{j,n}\in \R, (j,n)\in\N_m\ts \Z$.
Then the operator $\cJ=\os_1^m \cJ_j$, where $\cJ_j$ is a scalar
Jacobi operator, acting in $\ell^2(\Z)$ and is given by
$$
(\cJ_j y)_n=a_{j,n}y_{n+1}+a_{j,n-1}y_{n-1}+b_{j,n}y_n,\qq n\in\Z,\qq
  y=(y_n)_{n\in \Z}\in \ell^2(\Z).
$$
In this case each $\D_j$ is the standard Lyapunov function for
$\cJ_j $, and the properties of $\D_j$ is well known, see Theorem
\ref{16001} and \cite{T}. Thus $\r=1$, since $\o_j=\{j\}$.\BBox

Below we show the following identity (see \er{1020})
\[
 \lb{8000}
 D(z,\t)=c{\t^m}\prod_{j=1}^{mp}(z-\l_j(\t)),\qq  z,\t\in\C,\ \t\not=0,\ \ c=(-1)^m\det A_p,
\]
where $\l_j(\t)$ are zeros of the polynomial $D(z,\t)=0$ for
fixed $\t$. Define the set  $\L(\t)=\{\l_{n}(\t),n\in\N_{mp} \}, \t\in\C$.  Then we obtain $ \s(\cJ)=\bigcup_{j=1}^{mp}\s_j$,
where $\s_j=\l_j(\S^1), \S^1=\{\t\in\C:\  |\t|=1\}$ is the bounded close set. Ler  $\N_m^0=\{0,...,m\}$ for $m\ge 0$.
We  describe minimum spectral data which determine all Lyapunov functions
$\D_j, j\in \N_m$ for some Jacobi operator.

\begin{theorem}\lb{T2}
\it  Let $\vk_k\in \R, j\in\N_m^0$ satisfy
$\cos\vk_n\not=\cos\vk_j$ for all $j\not=n$.

\no\it i) Let $\L_0=\L(e^{i\vk_0})$ and let $\L_j\ss\L(e^{i\vk_j}),
\# \L_j=(m-j)p+1,j\in \N_m$ for some Jacobi operator $\cJ$. Then the spectral data $\L_j, j\in \N_m^0$,
 determine the polynomial $D(\cdot,\cdot)$, all Lyapunov functions
 $\D_j,j\in \N_m$ and the
spectrum $\s(\cJ)$ (counted according to its multiplicity).

{{\no\it ii)}} Let $\wt\L_1\ss \L_1$ and let $\#\wt\L_1=\#\L_1-1$.
Then there exist infinitely many Jacobi operators, having the same spectral data $\L_0$, $\wt\L_1$, $\L_j$, $j=2,...m$, but different determinants.

\end{theorem}

{\bf Example.} Note that if $e^{ik\vk_j}=1$ for some $k\in \N$, then $\l_n(e^{i\vk_j}), n\in \N_{mp}$ are k-periodic eigenvalues.
Consider the case $m=2$. Let $\vk_0=0, \vk_1=\pi$.
Then $\L_0$ is a set of all periodic eigenvalues
and $\L_1$ is a set of some anti-periodic eigenvalues
such that $\# \L_0=2p, \# \L_1=p+1$.
Let the set $\L_2=\{\wt\l\}$, where $\wt\l\in \s(\cJ)$ is a some
four-periodic eigenvalue, i.e., $\vk_2={\pi\/2}$.
 Due to Theorem \ref{T2}, the spectral data
$\L_0, \L_1, \L_2$ determine the polynomial $\det (\mM_p(z)-\t I_m)$,
 the spectrum $\s(\cJ)$ and the Lyapunov functions
$\D_1, \D_2$.\BBox

\begin{theorem}\lb{7100}\it
\no i) Let $p\ge3$. Then each sum
$\sum_{j=1}^{mp}\l_j^s(\t), s\in\N_{p-1}$ does not depend on $\t\in \C$. In
particular, the following identities hold
\[
 \lb{2010}
 \sum_{n=1}^{mp}\l_n(\t)=\sum_{n=1}^p\Tr b_n,\qqq
 \sum_{n=1}^{mp}\l_n^2(\t)=\sum_{n=1}^p\Tr (b_n^2+2a_na_n^\top),
 \qq \  all \qq \t\in\C.
\]
\no ii) The following estimate is fulfilled:
\[
 \lb{7110}
 \sum_{n=1}^{pm} \l_n^2(\t)\ge2pm(\det A_p)^{2\/{pm}},
\]
where the identity holds true iff $b_n\ev0$, $a_na_n^{\top}=(\det
A_p)^{2\/pm}I_m$ for all $n\in\N_{p}$.

\no iii) Let numbers $\vk_j\in\R, j\in\N_m$ satisfy
$\cos\vk_n\not=\cos\vk_j$ for all $j\ne n$. Then the eigenvalues
$\l_{nm+k}(e^{i\vk_j})=2\cos \frac1p(\vk_j+2\pi n)$ for all
$(j,n,k)\in\N_m^0\ts \N_{p-1}^0\ts\N_p$  iff $b_n=0,
a_na_n^\top=I_m$ for all $n\in\N_{pm}$ and $\prod_{n=1}^pa_n=I_m$.
\end{theorem}

We describe a priori estimates

\begin{proposition} \lb{T3}
\it Let $\|\cJ\|_{\iy}=\max\{|a_n(j,k)|$, $|b_n(j,k)|\}$, where
the matrices $a_n=\{a_n(j,k)\},\\ b_n=\{b_n(j,k)\}$. Then
the following estimates hold
\[
\lb{T3-1}
 \|\cJ\|_{\iy}\le\max\{|\l_0^+|,|\l_N^-|\}=
 \|\cJ\|\le(4m-1)\|\cJ\|_{\iy},
\]
\[
\lb{T3-2}
 \|\cJ\|_{\iy}+{|\l_{N}^-+\l_0^+|\/2}\le {\l_N^--\l_0^+\/2}\le (4m-1)\|\cJ\|_{\iy},
 \qq if \qq \sum_{j=1}^p\Tr b_j=0.
\]
\end{proposition}

\no {\bf Remark.} A priori estimates for scalar Jacobi operators
were obtained in \cite{BGGK}, \cite{K}, \cite{KKr}.

The plan of our paper is as follows. In Sect. 2 we  prove Theorem \ref{T1}. In Sect. 3 we prove Theorem \ref{T2}, \ref{7100}
and Proposition \ref{T3}. In Section 4 we consider examples
for the case $m=p=2$, where we construct the complex and real resonances.
 In Sect. 5 we shortly recall the well known results about
the properties of point spectrum and the absence of singular continuous spectrum. In the proof we use \cite{CK}, \cite{BBK}.

\section {The Lyapunov functions}
\setcounter{equation}{0}

We recall well-known results.
For any solution $y=(y_n)_{n\in\Z}$ of the equation
$a_n y_{n+1}+b_n y_n+a_{n-1}^\top y_{n-1}=z y_n$ we define
\[
 \lb{900}
 f_n=\left(\begin{array}{c} y_n \\ y_{n+1} \end{array}\right),\ \ and
 \qqq
 \cT_n=\ma 0 & I_m \\ -a_n^{-1}a_{n-1}^\top & a_n^{-1}(z-b_n) \am,\ \
 n\in\Z.
\]
Then $f_n$ satisfies $f_{n}=\cT_{n}f_{n-1}$. Thus the matrix-valued function
$\mM_n=\ma \vt_{n} & \vp_{n} \\
           \vt_{n+1} & \vp_{n+1}\am
$ satisfies $\mM_n=\cT_{n}\mM_{n-1}, \mM_0=I_{2m}$.
This gives the monodromy matrix $\mM_n=\prod_{j=1}^n \cT_j$, $n\ge1$,
where $\prod_{j=1}^{n}X_j=X_{n}...X_{1}$ for matrices $X_j$. We rewrite $\cT_n$ in the form
\[
 \lb{901}
 \cT_n=\ma I_m & 0 \\ 0 & a_n^{-1} \am\ma 0 & I_m \\ -I_m & z-b_n  \am
 \ma a_{n-1}^\top & 0 \\ 0 & I_m  \am= P_n^{-1}R_n\cT_n^0P_{n-1},
\]
where
\[
 \lb{902}
 P_n= a_{n}^\top \os I_m  ,\qq
 R_n=P_n\ma I_m & 0 \\ 0 & a_n^{-1} \am=a_n^\top \os a_n^{-1},
 \qq \cT_n^0=\ma 0 & I_m \\ -I_m & z-b_n  \am.
\]
Then using \er{901}, we get
\[
 \lb{903}
\mM_n=\prod_{j=1}^n P_j^{-1}R_j\cT_j^0P_{j-1}=P_n^{-1}\lt(\prod_{j=1}^nR_j\cT_j^0\rt)P_0=P_n^{-1}M_nP_0,\
 \ M_n=\prod_{j=1}^nR_j\cT_j^0.
\]
Matrices $R_j,\cT_j$, $j\ge1$ satisfy
\[
 \lb{904}
 R_j^\top J R_j=J,\ \ \ (\cT_j^0)^\top J \cT_j^0=J,\ \ \  J=\ma 0 & I_m \\ -I_m & 0  \am,
\]
which yields
\[
 \lb{905}
 M_n^\top J M_n=J,\ \ \ M_n^{-1}=-JM_n^\top J,\qq n\in \Z.
\]
Using \er{903}  we obtain
\[
M_n=P_n\mM_n P_0^{-1}
=\ma a_n^\top\vt_n(a_0^\top)^{-1}&a_n^\top\vp_n\\
\vt_{n+1}(a_0^\top)^{-1}&\vp_{n+1}\am,\qq
M_n^\top=\ma a_0^{-1}\vt_n^\top a_n&a_0^{-1}\vt_{n+1}^\top\\
\vp_n^\top a_n&\vp_{n+1}^\top\am
\]
and
\[
M_n^{-1}=-JM_n^\top J= \ma \vp_{n+1}^\top& -\vp_n^\top a_n\\
-a_0^{-1}\vt_{n+1}^\top & a_0^{-1}\vt_{n}^\top a_n\am .
\]
Due to $M_nM_n^{-1}=I_{2m}$ and $ M_n^{-1}M_n=I_{2m}$ we deduce that
\[
a_n^\top\vt_n (a_0^\top)^{-1}\vp_{n+1}^\top -a_n^\top\vp_n (a_0^\top)^{-1}\vt_{n+1}^\top =I_{m},\qqq
\vt_n (a_0^\top)^{-1}\vp_{n}^\top=\vp_n (a_0^\top)^{-1}\vt_{n}^\top,
\]
\[
\vt_n^\top a_n \vp_{n+1}-\vt_{n+1}^\top a_n \vp_{n}=a_0,\qq
\vp_{n+1}^\top a_n^\top \vp_{n}=\vp_n^\top a_n^\top\vp_{n+1},\qq
\vt_{n+1}^\top a_n^\top \vt_{n}=\vt_n^\top a_n^\top\vt_{n+1}.
\]
Recall the asymptotics of fundamental solutions
\[
\vp_{n+1}=z^{n}A_n+O(\l^{n-1}),\ \
 \vt_{n+1}=O(z^{n-1}), \ n\ge 1,\qq A_n=(a_1a_2...a_n)^{-1}
\]
as $z\to \iy$. Substituting last asymptotics into $\mM_p$ we obtain
\[
\lb{600} \mM_p(z)=z^p \ma 0&0\\ 0&A_p \am  +O(z^{p-1}).
\]
It is well known that the  determinant $D(z,\t)$ satisfies
\[
 \lb{009}
 D(z,\t)=\det(M(z)-\t I)=\t^{2m}+\sum_{j=1}^{2m}\t^{2m-j}\x_j(z)=
 \prod_{j=1}^{2m}(\t-\t_j(z)),\ \ z,\t\in\C,
\]
\[
\lb{15}
\x_0=1,\ \ \ \x_1=-T_1,\ \  \x_2=-{1\/2}(T_2+T_1\x_1),\ \ ... ,
\x_s=-{1\/s}\sum_0^{s-1}T_{s-j}\x_j,..,
\]
see p.331-333 [RS].
The identity \er{701} gives $\x_j=\x_{2m-j}$ for all $j\in \N_m$.
Recall that the Chebyshev polynomials $\mT_n, n\ge 1$ satisfy:
\[
\lb{CP}
\mT_n(\n)={\t^{n}+\t^{-n}\/2}=
2^{m-1}\sum_0^{[{n\/2}]}c_{n,j}\n^{n-2j},\ \ \
c_{n,j}=(-1)^j n{(n-j-1)!\/(n-2j)!j!}2^{n-2j-m} ,\ \
\]
$\n={\t+\t^{-1}\/2}$,  see [AS].
Then  \er{009} and the identity ${\t^{n}+\t^{-n}\/2}=\mT_n(\n)$ yield
\[
\lb{DCP}
{D(\t,z)\/(2\t)^m}=
{(\t^m+\t^{-m})\/2^m}+\x_1{(\t^{m-1}+\t^{1-m})\/2^m}+...+2^{-m}\x_{m}=\sum_{j=0}^m {\x_{j}\/2^{m-1}}\mT_j(\n),
\]
\[
\lb{Trn}
{1\/2}\Tr M^n(z)=\sum_1^m{\t_j^n+\t_j^{-n}\/2}=\sum_{j=1}^m\mT_n(\D_j(z)),\ \ \ z\in \C.
\]
The substitution of \er{CP} into the equality \er{DCP} gives
\[
\lb{14}
\F(\n,z)={D(\t,z)\/(2\t)^m}=\sum_0^m \f_j(z)\n^{m-j},
\]
\[
\lb{111}
\f_0=1,\ \ \f_1=c_{m-1,0}\x_1={\x_1\/2},\ \ \f_2=
c_{m,1}+c_{m-2,0}\x_2,\ \ \f_3=c_{m-1,1}\x_1+c_{m-3,0}\x_3,...,
\]
\[
\lb{112}
\f_{2n}=c_{m,n}+c_{m-2,n-1}\x_2+c_{m-4,n-2}\x_4+...+
c_{m-2n,0}\x_{2n},
\]
\[
\lb{113}
\f_{2n+1}=c_{m-1,n}\x_1+c_{m-3,n-1}\x_3+c_{m-5,n-2}\x_5+...+
c_{m-2n-1,0}\x_{2n+1}.
\]
Let  $a_n=I_m, b_n=0$ for all $n\in \Z$ and denote the corresponding Jacobi operator by $\cJ^0$. If $D^0$
is the corresponding determinant, then we obtain
% and the polynomial $\F^0=D^0(z,\t)/(2\t)^m$. Then  we have
$D^0(\t,z)=(\t^2+1-2\t\mT_p({z\/2}))^m$ and $\F^0={D^0(\t,z)\/(2\t)^m}=(\n-\mT_p({z\/2}))^m=\sum_0^mC_j^m(-\mT_p({z\/2}))^j\n^{m-j},
$
where $C_m^{N}={N!\/(N-m)!m!}$.

\no {\bf Proof of Theorem  \ref{T1}.}
i) We have the identity $
\F(\n,z)={D(\t,z)\/(2\t)^m}=\sum_0^m \f_j(z)\n^{m-j}$, see \er{14},
where the polynomials $\f_j$ are given by \er{111}-\er{113}.
$\D_1(z),..,\D_m(z)$ are
the roots of $\F(z,\n)=0$ for fixed $z\in \C$.
Then the statement i) follows from the well-known about the zeros
of a polynomial $\sum_0^m\n^{m-j}w_j(z)$, where $w_j$ is a polynomials
in $z$, see \cite{Fo}.

ii) We have $\D_j'(z)={1\/2}(1-\t^{-2}(z))\t'(z)\ne 0, z\in Y$,
since by the Lyapunov-Poincar\'e Theorem, $\t'(z)\ne 0$ for all $z\in Y$.

iii) Recall  that the resultant for polynomials $
f=\t^n+\a_{1}\t^{n-1}+..+\a_n, \  g=b_0\t^s+\b_{1}\t^{s-1}+..+\b_s $
is given by
\[
\lb{res} R(f,g)=\det \left(\begin{array}{cccccccc}
1& \a_1&..&\a_{n}&0&0&..&0\\
0&1& \a_1&..&\a_{n}&0&..&0\\
. &.&.&.&. &.&.&. \\
0&...&0&1& \a_1&&..&\a_{n}\\
\b_0&\b_1&..&\b_s&0&..&0\\
0&\b_0&\b_1&..&\b_{s}&0&..&0\\
. &.&.&.&. &.&.&. \\
0&...&0&\b_0&\b_1&&..&\b_s\\
\end{array}\right)
\begin{array}{c}
\left.\phantom{\begin{array}{c}
~\\~\\~\\~\end{array}}\right\} \mbox{ $n$ lines}\\
\left.\phantom{\begin{array}{c}
~\\~\\~\\~\end{array}}\right\}\mbox{ $s$ lines}
\end{array}.
\]
The discriminant of the polynomial $f$ with zeros $\t_1,..,\t_n$
is given by
$${\rm Dis} f=\prod_{i<j} (\t_i-\t_j)^2=(-1)^{{n(n-1)\/2}}R(f,f').
$$
Thus we have ${\rm Dis} \F_j(\l,\t)=\!\!\prod_{i<s, i,s\in
\o_j}\!\! (\D_i(\l)-\D_s(\l))^2=(-1)^{N_j(N_j-1)\/2}R(\F_j,\F_j')
$ is polynomial, since $\F_j(\l,\t)$ is the polynomial. Then the
function $\r=\prod_{1}^{N_0}{\rm Dis} \ \F_j$ is polynomial.

iv) Each gap has the form  $\g_n=(\l_n^-,\l_n^+)=\cap_{j=1}^m \g_{n,j}$, where
$\g_{n,j}=(\l_{n,j}^-,\l_{n,j}^+)\ss\R$ is some finite interval such that
$\D_j(z)\notin [-1,1]$ for  all $z\in \g_{n,j}$.
Note that $\D_j(\l_{n,j}^-)=\pm $ or $ \D_j(\l_{n,j}^+)=\pm 1$ or $\l_{n,j}^\pm$ is the branch
point of $\D_j(z)$, otherwise we have a contradiction.

v) Using \er{600} we define $F$ by
$$
 F(t)=t^p\mM_p=F_0+tF_1(t),\qqq
F_0=0_m\os A_p ,\ \ t=z^{-1},\qq
$$
where $F_1$ is a some matrix polynomial and $F_1(t)=\const +O(t)$
as $t\to 0$. Let $\wt\t_j(t), j\in\N_{2m}$ be the eigenvalues
of $F(t)$.
Recall that $\t_j^0, j\in \N_m$ are  eigenvalues of $A_p$.
Then, $\wt\t_j(t)$ are given by  convergent Puiseux series
(the Taylor series in $\t^{1/m_j}$)
$$
 \wt\t_j(t)=\t_j^0+\sum_{k=1}^\iy\a_{j,k}t^{{k\/m_j}},\
 \ |t|<r,\qq \sum_1^k m_j=m,  \qq  m_j\in \N,\ j\in w_s, s\in \N_k
$$
for some $r>0$, see p=4, \cite{RS}, where $w_1=\N_{m_1}, w_2=\N_{m_1+m_2}\sm w_1, ..., w_k=\N_m\sm \N_{m-m_k}$
and $\t_j^0=\t_{j'}^0$ for all $j,j'\in w_s$.
Here $ \wt\t_j(t), j\in w_s$ are the branches of one or more multivalued analytic functions with at worst algebraic branch points at $\t=0$
and satisfy
$$
%\[ \lb{910}
 \t_j(z)=z^p\wt\t_j(t)=
 z^p(\t_j^0+O(z^{-{1\/m_s}}))=\t_j^0z^p+O(z^{p-{1\/m_s}}),\
 \qq j\in w_s.
$$
 If each eigenvalue $\t_j^0$
is simple, then each $\wt\t_j(t),j\in\N_m$ is analytic for small $t$,
i.e.,  $m_s=1$. Thus we obtain $\t_j(z)=\t_j^0z^p+O(z^{p-1}), j\in\N_m$, which yields $2\D_j=\t_j(z)+\frac1{\t_j(z)}=\t_j^0z^p+O(z^{p-1})$ and
using \er{T1-2}, we get asymptotics for $\r_s$.

\section {Inverse Problems}
\setcounter{equation}{0}

In order to prove Theorem \ref{T2} we need the  Lemma 3.1-3.3.

\begin{lemma}\lb{11000}
Let $W=\{W_{r,j}\}_{j,r=0}^s$ be a $s\ts s$ matrix  with components $W_{r,j}=\cos j\vk_r$ for some $\vk_r\in \C, s\ge 1$, where
 $\cos\vk_r\ne\cos\vk_{r'}$ for all $r\not=r'$. Then
$\det W\ne 0$.
\end{lemma}
{\no\bf Proof.}
Using $\cos jz=\mT_j(\cos z)$, where the Chebyshev
polynomials  $\mT_j$ are given by \er{CP}, we obtain $W_{r,j}=\cos j\vk_r=\mT_{j}(\cos\vk_r)$.
Let $W y=0$ for some vector $y=(y_n)_0^s \in\C^{s+1},
y\ne 0$. Then we deduce that a polynomial $P(x)=\sum_{j=0}^s y_j\mT_{j}(x), P\ne 0$ has distinct zeros $\cos\vk_j, j\in \N_s^0$, which
gives a contradiction, since $\deg P\le\deg \mT_{s}=s$.
 Thus $W$ is invertible.\BBox

\begin{lemma} \lb{11001}
Let $c=(-1)^m\det A_p$. The polynomials $\x_j(z), j\in\N_{m}$, given by \er{009}, satisfy
\[
\x_j(z)=\x_{2m-j}(z)=O(z^{pj}),\qqq\
 \ \x_m(z)=cz^{pm}+O(z^{pm-1})\qq
as\ z\to \iy,\ ,
\]
\[
 \lb{1020}
 c\prod_{n=1}^{mp}(z-\l_n(\t))=
  {D(z,\t)\/\t^{m}}=(\t^m+\t^{-m})+\x_m(z)+\sum_{j=1}^{m-1}
 (\t^{m-j}+\t^{j-m})\x_j(z).
\]
\end{lemma}
{\no\bf Proof.} Using \er{701}, \er{009} and Viete's formulas, we
get
\[
 \lb{911}
 \x_j(z)=\x_{2m-j}(z)=(-1)^j\sum_{1\le
 i_1<..<i_j\le2m}\prod_{s=1}^j\t_{i_s}(z),\ \ 1\le j\le m.
\]
Asymptotics \er{T1-3} give $\t_j(z)=\t_j^0z^p+o(z^p)$,
$j\in\N_m$ and the Lyapunov-Poincar\'e Theorem yields
$\t_{m+j}(z)=\frac1{\t_j(z)}=O(z^{-p})$ as
$z\to\iy$, since each $\t_j^0\neq 0$. Then, using \er{911}, we
obtain
$$
\x_j(z)=O(z^{pj}),\qq
 \x_m(z)=(-1)^m\prod_{j=1}^m\t_j^0z^{pm}+o(z^{pm})=(-1)^m\det A_p
 z^{pm}+O(z^{pm-1}).
$$
Identity \er{1020} follows from last identity and \er{701},
\er{009}. \BBox

We need the following Lemma about polynomials.

\begin{lemma} \lb{11010}
Let $(z_n)_{1}^k\in \C^{k}$
and $(c_j)_{0}^s\in \C^{s+1}$ for some integers $s,k\ge 1$. Then there
exists a unique polynomial $r(z)$ satisfying for some  polynomial $g$
the following relations
\[
 \lb{a01}
 r(z)=z^k\sum_{j=0}^sc_jz^j+O(z^{k-1})\ \ {\rm\ as} \ z\to \iy\qq {\rm\ and\ }\
 r(z)=h(z)g(z),\ \ h(z)=\prod_{n=1}^k(z-z_n),
\]
\end{lemma}
{\no\bf Proof.} Introduce the linear space of polynomials
$
 \cP_s=\{g:\ \deg g\le s\},\ \ \dim\cP_s=s+1$.
 Note that $g\in \cP_s$.
Define the linear operator $\cA:\cP_{s}\to\C^{s+1}$ by
$$
 (\cA   g)_j={1\/n!} {d^n\/dz^n}\rt(h(z) g(z)\rt)|_{z=0},\ \
 n=k+j-1,\  j\in \N_{s+1},\ \ \ g\in\cP_{s}.
$$
We rewrite $\cA   g$ in the form $(\cA  g)_j=c_{j-1},\ \ j\in
\N_{s+1}$, where $h(z)g(z)=z^k\sum_{j=0}^sc_jz^j+O(z^{k-1})$ as
$z\to \iy$. We deduce that $\cA:\cP_{s}\to\C^{s+1}$ is an
isomorphism, since $\dim \cP_{s}=\dim \C^{s+1}=s+1$ and if $\cA
g=0$, then $g=0$. If $g=\cA^{-1}(c_n)_{0}^s$ and $r(z)=h(z)g(z)$,
then $ r(z)=z^k\sum_{j=0}^sc_jz^j+O(z^{k-1})$ as $z\to\iy$.
Moreover, the polynomial $r$ is unique, since $\cA$ is an
isomorphism. \BBox

\no{\bf Proof of Theorem \ref{T2}.} Recall that $c=(-1)^m\det A_p$. Define
the polynomials $\zeta_0(z)\ev c^{-1}$, $\zeta_j(z)\ev
c^{-1}\x_j(z)$, $j\in\N_m$. Lemma \ref{11001} yields $\deg\zeta_j\le pj$ and then $\zeta_j(z)=\sum_{n=0}^{pj}\zeta_{j,n}z^n$ for some
$\zeta_{j,n}\in\R$. Introduce the sets
$
 K_s=\{p(m-s-1)+1,...,p(m-s)\},\ \
 s\in \N_{m-1}^0,\ \ K_{m}=\{0\},
$
which satisfy
\[
\lb{dks}
\cup_{s=j}^m K_s=\N_{p(m-j)}^0,\qq   K_r\cap K_s=\es, for\ \ all \
\ r\ne s.
\]
 Using \er{1020} and \er{dks} and setting $h_j=\t^{-j}+\t^j$, we obtain
\[
 \lb{7000}
q(z,\t)={D(z,\t)\/c\t^{m}}=\prod_{n=1}^{pm}(z-\l_n(\t))=
 \sum_{j=0}^m h_{m-j}\z_j(z) =\sum_{j=0}^mh_j\z_{m-j}(z)
\]
$$
=\sum_{j=0}^mh_j\sum_{n=0}^{p(m-j)}\z_{m-j,n}z^n=
 \sum_{j=0}^m h_j\sum_{s=j}^m\sum_{n\in K_s}\z_{m-j,n}z^n
 =\sum_{j=0}^m\sum_{s=j}^m\sum_{n\in K_s}h_j\z_{m-j,n}z^n
 $$
$$
 =\sum_{s=0}^m\sum_{j=0}^s\sum_{n\in  K_s}h_j\z_{m-j,n}z^n
 =\sum_{s=0}^m\sum_{n\in K_s}\sum_{j=0}^sh_j\z_{m-j,n}z^n=
 \sum_{s=0}^m\sum_{n\in K_s}\eta_n(\t)z^n=\sum_{n=0}^{pm}\e_n(\t)z^n,
$$
where
\[
 \lb{7001}
 \eta_{n}(\t)=\sum_{j=0}^sh_j\z_{m-j,n},\
 \ n\in K_s,\ \ s\in\N_m^0.
\]
Substituting $\t=e^{i\vk_r}, r\in\N_s^0$ into the
identity \er{7001} we get
\[
 \lb{7020}
\wt\e_n=2W_s y_n, \qq \wt\e_n=(\e_{n}(e^{i\vk_r}))_{r=0}^s,\
 y_n=(\z_{m-j,n})_{j=0}^s,\ \ W_s=(\cos j\vk_r)_{r,j=0}^s,
\]
for all $(n,s)\in K_s\ts\N_m^0$.
Due to Lemma \ref{11000}, each matrix $W_s, s\in\N_m^0$ is
invertible. Thus \er{7020}, \er{7001} give
\[
\lb{11002}
 \eta_n(\t)=<(\t^j+\t^{-j})_{j=0}^s,(2W_s)^{-1}\wt\e_{n}>,\qq
 (n,s)\in K_s\ts\N_m^0
\]
where $<\cdot,\cdot>$ is a scalar product in $\R^s$.
Then the functions $\e_n(\cdot), n\in N_s$ are determined by $\wt\e_n=(\e_{n}(e^{i\vk_r}))_{r=0}^s$.  We will use this fact
to determine all functions $\e_n(\cdot), n\in \N_m$.
In order to describe recovering we need the following simple fact.

\begin{lemma} \lb{a001}
Let $0\le k\le m-1$. Then numbers
$\e_n(e^{i\vk_j}), (n,j)\in \N_{pm}^0\ts\N_k^0$, the
functions $\e_n$, $n\in\cup_{j=0}^{k} K_j$, the set
$\L_{k+1}$ and the numbers $\vk_j$, $j\in\N_m^0$
determine the numbers
$\e_n(e^{i\vk_j}), (n,j)\in \N_{pm}^0\ts\N_{k+1}^0$,
and the functions $\e_n$, $n\in\cup_{j=0}^{k+1} K_j$.
\end{lemma}
{\no\bf Proof.} Let $S_{k+1}=(\e_n(e^{i\vk_{k+1}}))_{n\in \cup_{j=0}^{k}
 K_j}=(\e_n(e^{i\vk_{k+1}}))_{n=p(m-k-1)+1}^{pm}$.
%q(z,\t)\ev c^{-1}\t^{-m}D(z,\t),\qqq
Using \er{7000}, we see that the elements of $\L_{k+1}$
are some zeros of the polynomial $q(z,e^{i\vk_{k+1}})$ and
recall that $\# \L_{k+1}=p(m-k-1)+1$. Also, using \er{7000}, we see that
the components of vector $S_{k+1}$ are the coefficients of the
polynomial $q(z,e^{i\vk_{k+1}})$. Then, due to Lemma \ref{11010},
 $\L_{k+1}$ and $S_{k+1}$ uniquely determine the
polynomial $q(z,e^{i\vk_{k+1}})$ and, using \er{7000}, we
determine $(\eta_n(e^{i\varkappa_{k+1}}))_{n=0}^{pm}$.
Substituting obtained values $(\eta_n(e^{i\varkappa_j}))_{n\in
K_{k+1}}$, $j\in\N_{k+1}^0$ into the \er{11002}, we
determine $\eta_n(\t)$, $n\in K_{k+1}$. \BBox

Now we describe the recovering procedure. Recall that $q(z,\t)= {D(z,\t)\/c\t^m}$. Substituting the
elements of the set $\L_0$ into the \er{7000}, we determine
$q(z,e^{i\vk_0})$ and
$(\eta_n(e^{i\vk_0}))_{n=0}^{pm}$. Then \er{11002}
determines all functions $\e_n(\t), n\in K_0$. Step by step, using Lemma
\ref{a001} for $k\in\N_{m-1}^0$, we determine $\eta_n(\t)$,
$n\in\cup_{j=0}^mK_j=\N_{pm}^0$. Then \er{7000} gives $q(z,\t)$. Using \er{009} we get $cq(z,\t)=D(z,\t)=\t^{2m}+O(\t^{2m-1})$ as $\t\to\iy$,  then the polynomial $q(z,\t)$ uniquely determine constant $c$ and
$D(z,\t)$.

ii) Assume that $a_n, b_n, n\in\Z$ are diagonal,
and each $a_n>0, n\in\Z$. Then $\cJ=\os_1^m \cJ_j$, where $\cJ_j$ are scalar Jacobi operators. Let $\D_j$ be the Lyapunov functions
and let $\wt D_j$ be the corresponding determinant for
the operators $\cJ_j, j\in\N_m$. Then we deduce that
\[
 \lb{20000}
 D(z,\t)=\prod_{j=1}^m \wt D_j,\
\wt D_j (z,\t)=\t(\t+\t^{-1}-2\D_j(z)).
\]
By Theorem \ref{16001}, the Lyapunov function $\D_m(z)$ satisfies $ \D_m'(z_j^0)=0$ and $(-1)^{p-j}\D_m(z_j^0)=|\D_m(z_j^0)|\ge1,\ \ j\in \N_{p-1}$ for some $z_1^0<...<z_p^0$. Then the polynomial $q_s(z)=s(\D_m(z)-\cos\vk_0)+\cos\vk_0, s\ge 1$ satisfies
$$
(-1)^{p-j} q_s(z_j^0)= s\lt(|\D_m(z_j^0)|+(-1)^{p-j}\lt(-1+{1\/s}\rt)\cos\vk_0\rt)\ge 1,\qq (q_s)'(z_j^0)=0,\ \ j\in \N_{p-1}.
$$
 Using $\D_m=c_mz^p+O(z^{p-1}), c_m>0$ as $z\to \iy$, we get $q_s(z)=c_msz^{p-1}+O(z^{p-1})$, $c_m>0$. Due to Theorem  \ref{16001},
$q_s$ is a Lyapunov function for some scalar $p$-periodic Jacobi operator
$\cJ_m^s$. Consider the $p$-periodic $m$-dimensional Jacobi
operator $\cJ^s=\cJ_m^{s}\os\os_1^{m-1}\cJ_j$. Then the
determinant for $\cJ^s$ is given by
\[
 \lb{20001}
D_s(z,\t)=\t^{m}(\t+\t^{-1}-2q_s(z))
\prod_{ j=1}^{m-1}(\t+\t^{-1}-2\D_j(z)).
\]
Using \er{20000}, \er{20001}, we see that $\L_0(\cJ)=\L_0(\cJ^s)$, where $\L_j(\cJ),j=0,1,2..$ is  the set $\L_j$ for $\cJ$.  Each Lyapunov function $\D_{j_1}(z), j_1\in\N_{m-1}$ satisfies
$$
 \D_{j_1}(z)-\cos\vk_{j}=c_{j_1}\prod_{n=1}^p(z-z_{n+pj_1,j}),\ \
 \qq c_{j_1}>0,\ j\in \N_m,\ \ for\ some \ \ z_{n+pj_1,j}\in \R.
$$
 Then using \er{20000}, \er{20001},
we deduce that $z_{n+p,j}, n\in \N_{p(m-1)}$ are zeros of each
polynomial $D_s(z,e^{i\vk_j}), j\in\N_m$, $s\ge1$, where $
D_1(z,\t)=D(z,\t)$. We take $\L_j(\cJ^s)=\L_j(\cJ)=\{z_{p+n}, n\in\N_{p(m-j)+1}\}$, $j=2,..,m$ and
$\wt\L_1(\cJ^s)=\wt\L_1(\cJ)=\{z_{n+p,1}, n\in \N_{p(m-1)}\}$,
where $\wt\L_1(\cJ)$ is  the set $\wt\L_1$ for $\cJ$. But we have
$D_s\ne D_{1}=D$ for any $s>1$, since $q_{1}=\D_m\ne q_s$.
  \BBox

\no {\bf Proof of Theorem \ref{7100}.} i) Let $g_n(\t)=(-1)^n
\sum_{1\le j_1<..<j_n\le mp}\prod_{s=1}^n\l_{j_s}(\t)$. Lemma
\ref{11001} implies $\x_j(z)=O(z^{pj})$. Then using \er{1020} for
fixed $\t\neq\t_1=1$, $\t\ne 0$, we obtain
$$
 c\sum_{j=1}^{mp}z^{mp-j}(g_j(\t_1)-g_j(\t))=
 c\prod_{j=1}^{mp}(z-\l_j(\t_1))-c\prod_{j=1}^{mp}(z-\l_j(\t))=
$$
$$
 =(\t_1^m+\t_1^{-m}-\t^m-\t^{-m})+
 \sum_{n=1}^{m-1}
 (\t_1^{m-n}+\t_1^{n-m}-\t^{m-n}-\t^{n-m})\x_n(z)=O(z^{(m-1)p}),
$$
which yields $g_j(\t_1)=g_j(\t)$ for all $j\in\N_{p-1}$, since
$c\not=0$. The Newton formulas give
$$
 u_1=-g_1,\ \ \ u_n=-n g_n-\sum_{j=1}^{n-1}g_ju_{n-j},\qq
  where\qq u_n(\t)=\sum_{j=1}^{mp}\l_j^n(\t),
$$
which yields that $u_n(\t)=u_n(\t_1)=u_n(1)$ for all
$(n,\t)\in\N_{p-1}\ts\C$. Due to \er{3001}, we get
$$
\sum_{j=1}^{mp}\l_j^s(\t)=\Tr L(\t)^s, s\ge 1,\ \qq
 \sum_{j=1}^{mp}\l_j(\t)=\Tr L(\t)=\sum_{j=1}^N\Tr b_j
$$
\[
\lb{t11}
\sum_{j=1}^{mp}\l_j^2(\t)=\Tr L(\t)^2=\Tr
L(\t)L^*(\t)=2\sum_{j=1}^p \|a_j\|^2+\sum_{j=1}^p\|b_j\|^2.
\]

ii) Let $\ve_{j,n}$, $j\in\N_m$ be eigenvalues of the
matrix $a_na_n^\top, n\in\N_p$. Then \er{t11} gives
\[
 \lb{7101}
 \sum_{n=1}^{pm} \l_n^2(\t)\ge\sum_{n=1}^{p}\Tr2a_na_n^\top=
 2\sum_{n=1}^{pm}\sum_{j=1}^m\ve_{j,n}\ge
2pm  \lt(\prod_{j,n}\ve_{j,n}\rt)^{1\/pm}=2pm(\det
 A_p^2)^{1\/pm}.
\]
Recall that the identity $\sum_{n=1}^{p}\sum_{j=1}^m\ve_{j,n}=
pm  \lt(\prod_{j,n}\ve_{j,n}\rt)^{1\/pm}$ holds true
iff each $\ve_{j,n}^{pm}=\det A_p^2$,\ $(j,n)\in \N_m\ts \N_p$.
Then the statement ii) holds true.

iii) {\bf Sufficiency}.
Let $a^0=(a_n^0)_{n\in \Z}, b^0=(b_n^0)_{n\in \Z},$
where $a_n^0=I_m, b_n^0=0$ for all $n\in \Z$.
Let the Jacobi operator $\cJ=\cJ(a,b^0)$, where
each matrix $a_n, n\in \N_p$ is unitary  and $\prod_1^p a_n=I_m$. Define matrices $c_n$ by $c_1=I_m$,
$c_{n+1}=a_n^\top c_n$, $n\in\Z$. Thus $c_{n+p}=c_n$ for all $n\in\Z$ and
$c_n$ are unitary matrices, since $a_n$ are unitary. Then
$U=\diag(c_n)_{n\in\Z}$ is unitary operator in $\ell^2(\Z)^m$ and
$U\cJ(a,b^0) U^{-1}=\cJ(a^0,b^0)=\cJ^0$, which yields
$$
 \det (M_p-\t I_{2m})=\t^{m}\prod_{j=1}^m(\t+\t^{-1}-2\mT_p(z/2))=
 \det (M_p^0-\t I_{2m}),
$$
where $M_p, M_p^0$ are corresponding determinants and $\mT_p$ is the Chebyshev polynomial, see \er{CP}.  Using
\er{8000} and $\mT_p(z)=\cos(p\arccos z)$, we obtain
$\l_{nm+k}(e^{i\vk_j})=2\cos \frac1p(\vk_j+2\pi n)$ for all
$(j,n,k)\in\N_m^0\ts \N_{p-1}^0\ts\N_p$. Also, using i), we obtain
\[
 \lb{a0001}
 \sum_{n=1}^{pm}\l_n(e^{i\vk_j})^2=\sum_{n=1}^p\Tr ((b_n^0)^2+2a_n^0{a_n^0}^\top)= 2pm,\ \ j\in\N^0_m.
\]

{\bf Necessity}. Let $\l_{nm+k}(e^{i\vk_j})=2\cos \frac1p(\vk_j+2\pi n)$
for all $(j,n,k)\in\N_m^0\ts \N_{p-1}^0\ts\N_p$. Then Theorem
\ref{T2}.ii  uniquely determines
\[
 \lb{a0002}
 \t^{-m}D(z,\t)=\prod_{j=1}^m(\t+\t^{-1}-2\mT_p(z/2))
\]
(see above {\bf Sufficiency}). Identities \er{8000} give $\det
A_p=1$. Using \er{a0001} and Corollary \ref{7100}.ii, we obtain
$b_n=0$, $a_na_n^{\top}=I_m$. Due to \er{T1-1}, \er{T1-3} and
\er{a0002}, we obtain
$$
 \D_j(z)=\frac12(\t_j(z)+\t_j^{-1}(z))={\t_j^0z^p\/2}+
 O(z^{p-1}),\ \ \D_j(z)=\mT_p(z/2)={z^p\/2}+O(z^{p-1}),
$$
then $\t_j^0=1$, $j\in\N_m$, i.e.,  all eigenvalues of
matrix $A_p$ are equal $1$ and all eigenvalues of matrix
$A_p^{-1}=a_1...a_p$ are equal $1$, then $A_p=I_m$, since $a_n$
($a_na_n^\top=I_m$) and $A_p$ are unitary.
\BBox

\no {\bf Proof of Proposition \ref{T3}.}
For the self-adjoint operator $\cJ$ we obtain $\|\cJ\|=\max\{|\l_0^+|,|\l_0^-|\}$. Moreover,
the operator $\cJ$ has the form $\cJ=\sum_{2m-1}^{1-2m}\cJ_j$,
where $\cJ_j$ is some diagonal matrix.
Note that $\cJ_j, \cJ_s, j\ne s$ have different diagonals. Using the identity
$\|\cJ_j\|=\|\cJ_j\|_{\iy}$ and the simple estimate $\|\cJ\|_{\iy}\le \|\cJ\|$ we obtain
$$
\|\cJ\|_{\iy}\le \|\cJ\|\le\sum_{-2m+1}^{2m-1}\|\cJ^{(j)}\|\le(4m-1)\|\cJ\|_{\iy},
$$
which yields \er{T3-1}. If $\sum_1^p\Tr b_j=0$, then $\sum\l_j(\t)=0$, which yields $\l_0^+\le 0\le \l_0^-$ and $\l_0^--\l_0^+=|\l_0^+|+|\l_{N_G}^-|$. Thus all these stimates gives \er{T3-2}.\BBox

\section {Examples for the case $p=2$, $m=2$}
\setcounter{equation}{0}

We consider the case $p=2$, $m=2$, where real $2\ts 2$ matrices $a_n, b_n$ satisfy
  $$
a_{n+2}=a_n=\ma 1 & \b_{2n+1} \\ 0 & 1\am,\ \ b_{n+2}=b_n=\ma \a_{2n} & \b_{2n} \\
 \b_{2n} & \a_{2n+1} \am,\ \ \a_{n+4}=\a_n,\b_{n+4}=\b_n.
$$
The determinant $D$, the function $\r$ and the Lyapunov functions are given by
$$
D(\cdot,\t)=\t^4-T_1\t^3+T\t^2-T_1\t+1=(\t^2-2\D_1\t+1)
(\t^2-2\D_2\t+1)=(2\t)^2\F(\cdot,\n),
$$
$$
T={T_1^2-T_2\/2},\qq
\r={-T_1^2+2T_2+8\/4},\qq
\D_1={T_1-\sqrt{4\r}\/4},\qq \D_2={T_1+\sqrt{4\r}\/4},
$$
$$
 T_1=2z^2+T_{11}z+T_{10},\qq T_2=z^4+T_{23}z^3+..+T_{20},\qq
 \F(\cdot,\n)=\n^2-{T_1\/2}\n+{T\/4}-{1\/2},
$$
where
\[
\lb{4m} T_{11}=-\sum_0^3\a_n+(\b_0+\b_2)\b_3+\b_1(\b_0+\b_2),
$$
$$
T_{10}=-4+ 2\b_0\b_2+2\b_1\b_3
+\prod_0^3
\b_n+\a_0\a_2+\a_1\a_3-\b_0(\a_2\b_3+\a_3\b_1)-\b_2(\a_0\b_1+\a_1\b_3),
\]
\[
\lb{4B} T_{23}={-\sum_0^3 \a_n},\qq
 T_{22}=\sum _{0\le j<k\le 3}\a_j\a_k-\sum_0^3\b_n^2- 4,
$$
$$
T_{21}=(\a_2+\a_3)\b_0^2+(\a_0+\a_3)\b_1^2+(\a_0+\a_1)\b_2^2+(\a_1+\a_2)\b_3^2
$$
$$
-2\b_0(\b_1+\b_3)-2\b_2(\b_1+\b_3)+
(\a_1+\a_3)(2-\a_0\a_2)+(\a_0+\a_2)(2-\a_1\a_3),\ \
$$
$$
T_{20}=6+\b_0^2\b_2^2+\b_1^2\b_3^2-4\b_1\b_3-4\b_0\b_2 + \prod_0^3
\a_n-\a_0(\a_1\b_2^2+\a_3\b_1^2)-\a_2(\a_1\b_3^2+\a_3\b_0^2)
$$
$$
-2\a_0\a_2- 2\a_1\a_3+ +2\b_0(\a_2\b_3+\a_3\b_1)
+2\b_2(\a_0\b_1+\a_1\b_3).
\]
Below we consider some specific cases.

\no {\bf 1. All $a_n,b_n$ are diagonal
matrices.} In this case all $\b_n=0$, the function $\r=0$ and
$$
\D_1=\frac12\lt(z-{\a_0+\a_2\/2}\rt)^2-\frac12\lt({\a_0-\a_2\/2}\rt)^2-1,
\ \
%$$$$
\D_2=\frac12\lt(z-{\a_1+\a_3\/2}\rt)^2-\frac12\lt({\a_1-\a_3\/2}\rt)^2-1.
$$

\no {\bf 2. The constant coefficients.} Let
$\a_n=0,\b_n=\b$ for all $n\in\Z$. Then
$$
T_1=2z^2+4\b^2z+\b^4+4\b^2-4,\ \ \ \r=\b^2(2z+\b^2)^2(4z+\b^2+8).
$$
The function $\r$ has 3 zeros $z_1,z_2,z_3$.
If $\b=0$, then $\r= 0, \D_1\ev\D_2={z^2\/2}$.

If $\b\ne 0,\pm \sqrt 8$, then $\r$ has  the zeros  $z_1=z_2=-{\b^2\/2}, z_3=-{\b^2+8\/4}$.

If $\b=\pm \sqrt 8$, then $\r$ has  the zeros $z_1=z_2=z_3-4$.

We have the following equations for periodic and anti-periodic eigenvalues
$$
4(\D_1-1)(\D_2-1)=(z+2)^2((z-2)^2-4\b^2)=0,\qq
4(\D_1+1)(\D_2+1)=(z^2-2\b^2)^2=0.
$$
Hence for each $\b\ne 0$ there exist two double anti-periodic
eigenvalues  $-\b\sqrt 2, \b\sqrt 2$;
 one double periodic eigenvalue $-2$ and two simple periodic eigenvalues
$2(1-b), 2(1+b)$. If $\b=\pm 2$, then there exists the periodic
eigenvalue $-2$ of multiplicity 3.

\no {\bf 3. Example $(\a_{0},\a_1,\a_2,\a_3)=(1,0,-1,0)$;
$(\b_0,\b_1,\b_2,\b_3)=(t,0,0,0)$, $t\in\R\sm\{0\}.$} In this case
we have
$$
 T_1=2z^2-5,\qq T_2=z^4-(t^2+5)z^2-t^2z+8,
 \qq
 4\r_t=4t^2z^2+4t^2z+1.
$$
The Lyapunov functions and resonances $r_t^\pm$ (zeros  of $\r$) are given by
$$
 \D_{1,t}={2z^2-5-\sqrt{4\r_t}\/4},\qq
 \D_{2,t}={2z^2-5+\sqrt{4\r_t}\/4},\qq
 r_t^\pm={-1\pm \sqrt{1-t^{-2}}\/2}.
$$
Each resonance $r_t^\pm, t\in(0,1)$ is complex.

Let $t=1$. Then  the Lyapunov functions are given by
$\D_{1,1}={z^2-z-3\/2}, \D_{2,1}={z^2+z-2\/2}$ and the spectral
bands $\s_t^{j}=\{z\in \R: \D_{j,t}(z)\in [-1,1]\}, j=1,2$ have the forms
$$
 \s_1^{1}=[\l^{1}_1,\l^1_2]\cup[\l^1_3,\l^1_4]=
 \lt[-\frac{\sqrt{21}-1}2,-\frac{\sqrt5-1}2\rt]\cup
 \lt[\frac{\sqrt{5}+1}2,\frac{\sqrt{21}+1}2\rt],
$$
$$
 \s_1^{2}=[\l^{2}_1,\l^2_2]\cup[\l^2_3,\l^2_4]
 =\lt[-\frac{\sqrt{17}+1}2,-1\rt]\cup\lt[0,\frac{\sqrt{17}+1}2\rt].
$$
where $\l^{j}_k, j=1,2, k=1,4$ are periodic eigenvalues,
$\l^{j}_k, j=1,2, k=2,3$ are   anti-periodic eigenvalues and the resonances  $r_1^\pm=-{1\/2}$ lies on the gap $(-1,0)$.

Let $t=1+\ve$ for some small $\ve>0$. Then  the spectral bands
have the forms
$$
 \s_t^j=[\l_{1,t}^j,\l^j_{2,t}]\cup[\l_{3,t}^j,\l_{4,t}^j],\ \
 \l_{i,t}^j=\l_i^j+o(1),\ \ (i,j)\in \N_4\ts \N_2,\
 r_t^\pm={1+o(1)\/2} \qq as \qq  t\downarrow1.
$$
 Note that the "resonance" gap $\g_t=(r_t^-,r_t^+)\ne \es$, since  $r_t^-<r_t^+$ and $\s(\cJ)\cap\ol{\g_t}=\es$ as
$t\downarrow1$. Thus the interval $\g_t\ss (\l_{2,t}^2,\l_{3,t}^2)$, where $(\l_{2,t}^2,\l_{3,t}^2)$
is a gap in the spectrum of $\cJ$. This gives that "resonance gap"
 arises in the gap and the end points of the spectral bands
 are periodic or anti-periodic eigenvalues.

 If $t=1$, then the resonances $r_1^+=r_1^-$ are not the branch points.
 If $t>1$, then the resonances $r_1^-<r_1^-$ are the real  branch points.
If $t\in (0,1)$, then the resonances $r_1^\pm \in \C_\pm$ are the complex branch points.

%Пусть $|t|>1$. Тогда функция $\r(\l)<0$, если $\l\in(t_1,t_2)$ и
%$\r\ge0$, если $\l\in\R\sm(t_1,t_2)$, тем самым интервал
%$(t_1,t_2)$ есть резонансная лакуна, т.е. участок в котором не
%может быть точек спектра оператора $\mJ$. Все резонансы
%вещественны, однократны и образуют края резонансной лакуны. Если
%же $|t|<1$, то функция $\r(\l)>0$ при всех $\l\in\R$, т.е. в
%данном случае резонансных лакун нет и все резонансы комплексные.
%Если $t=1$, то $\r$ есть точный квадрат линейной функции и в этом
%случае функции Ляпунова $\D_{1,2}$ являются полиномами и не имеют
%ветвлений в $\C$.
%
%{{\no}} Далее посмотрим на положение резонансной лакуны
%относительно спектра $\s(\mJ)$. Если $t^2\ge\frac{\sqrt2+1}2$, то
%левый край резонансной лакуны 'попадает' в спектр, т.е. для
%некоторого $\ve>0$ будет $[t_1-\ve,t_1]\ss\s(\mJ)$ и конечно
%$(t_1,t_2)\cap\s(\mJ)=\varnothing$, при этом правый край в спектр
%не попадает, т.е. $[t_2,t_2+\ve]\cap\s(\mJ)=\varnothing$ при
%некотором $\ve>0$. Далее, при $1<t^2<\frac{\sqrt2+1}2$ резонансная
%лакуна невырождена, но оба её края не являются точками спектра
%(она существует сама по себе). При $t^2=1$ она вырождается и далее
%при $t^2<1$ её 'нет на вещественной оси'.

\no {{\bf  Example 4}. Let $(\a_{0},\a_1,\a_2,\a_3)=(0,1,0,1)$;
$(\b_0,\b_1,\b_2,\b_3)=(t,0,0,0)$, $t\in\R$}. Then we obtain
$$
 T_1=2z^2-2z-3,\ \ T_2=z^4-2z^3-(3+t^2)z^2+(4+t^2)z+4,
 \ \
$$
$$
 4\r_t=(2z-1)^2+4t^2(z^2-z),\ \ r_t^\mp={1\/2}\mp {t\/2\sqrt{t^2+1}}\in \R.
$$
Here $r^\pm_t$ are the zeros (resonances) of $\r_t$ and
 if $t\ne 0$, then $r_t^-<r_t^+$. Thus we have the resonance gap $\g_t=(r_t^-,r_t^+)\ss (0,1), t\ne 0$. The Lyapunov functions are given by
$$
\D_{1,t}=\frac14\lt(2z^2-2z-3-\sqrt{4\r_t(z)}\rt), \qqq
\D_{2,t}={1\/4}\lt(2z^2-2z-3+\sqrt{4\r_t(z)}\rt).
$$

Let $t=0$. Then the Lyapunov functions have the forms
$ \D_{1,0}={z^2-2z-1\/2},\ \ \D_{2,0}={z^2-2\/2}$
 and $\D_{1,0}({1\/2})=\D_{2,0}({1\/2})=-{7\/8}$. The spectral bands $\s_t^{j}=\{z\in \R: \D_{j,t}(z)\in [-1,1]\}, j=1,2$ have the forms
$$
\s_0^1=[-1,3], \ \ \s_0^2=[-2,2].
$$
The end points of spectral bands $\s_0^1,\s_0^2$ are periodic  eigenvalues.

 If $t\in(0,\ve)$ for some small $\ve$, then   we obtain
$$
 \s_t^j=[\l_1^j(t),r^-_t]\cup[r^+_t,\l_2^{j}(t)],\ \  r^-_t<r^+_t,
 \qq r_t^\mp={1+o(1)\/2},
$$
$$
 \l_{1,t}^{1}=-1+o(1),\ \ \ \l_{2,t}^{1}=3+o(1),\ \ \l_{1,t}^{2}=-2+o(1),\  \ \ \l_{2,t}^{2}=2+o(1)\qq as \qq t\to0.
$$
This yields $\s(\cJ)=\s_t^1\cup \s_t^2=[\l_{1,t}^{2},r_t^-]\cup[r_t^+,\l_{4,t}^{1}]$, where   $(r^-_t,r^+_t)\ne \es$ is the  resonance gap, since $r^-_t<r^+_t$,
and $\l_{2,t}^{1}, \l_{2,t}^{1}$ are periodic eigenvalues.

\section {Appendix}
\setcounter{equation}{0}

Let $ \mH=L^2\lt([0,2\pi),{dx\/2\pi},\cH\rt)=
 \int_{[0,2\pi)}^{\os}\cH{dx\/2\pi}
$
be  a constant fiber direct integral, where $\cH=\C^{pm}$.
It is the Hilbert space of square integrable $\cH$-valued functions.
Let $\mS(\Z)$ be the set of all compactly supported functions $(f_n)_{n\in \Z}\in l^2(\Z)^m$.

\begin{lemma} \lb{Tdi}
The operator $U:l^2(\Z)^m\to\mH$, given by
\[
 \lb{5001}
 (Uf)_j(x)=\sum_{n=-\iy}^{\iy}e^{-ixn}f_{j+np},\qqq
  (x,j)\in [0,2\pi ]\ts \N_p,
\]
 is well defined for $\mS(\Z)$
and uniquely extendable to a unitary operator. Moreover,
\[
 \lb{5004}
 U\cJ U^{-1}=\int_{[0,2\pi)}^{\os}L(e^{ix}){dx\/2\pi},\qqq
 L(\t)=\left(\begin{array}{ccccc} b_1 & a_1 & 0 & ... & \t^{-1} a_p^\top \\
                                  a_1^\top & b_2 & a_2 & ... & 0 \\
                                  0 & a_2^\top & b_3 & ... & 0 \\
                                  ... & ... & ... & ... & ... \\
                                  \t a_p & 0 & ... & a_{p-1}^\top & b_p
                                  \end{array}\right),
\]
\[
 \lb{Tdi-3}
 \s(\cJ)=\s_{ac}(\cJ)=\bigcup_{x\in[0,2\pi)}\s(L(e^{ix}))=
 \bigcup_{x\in[0,2\pi)}\cup_{n=1}^{pm}\{\l_n(e^{ix})\}.
\]
\end{lemma}

\no{\bf Proof.}
We use  standard arguments from \cite{RS}.
For  $f\in \mS(\Z)$ the sum \er{5001} is
clearly convergent. For such functions $f$ we compute
$$
\|Uf\|^2=\int_0^{2\pi}\|(Uf)(x)\|^2{dx\/2\pi}=
 \int_0^{2\pi}\sum_{j=1}^p\lt(\sum_{n=-\iy}^{\iy}e^{-ixn}f_{j+np},
 \sum_{s=-\iy}^{\iy}e^{-ixs}f_{j+sp}\rt){dx\/2\pi}
$$
$$
 =\sum_{j=1}^p\lt(\sum_{n=-\iy}^\iy\sum_{s=-\iy}^\iy
 f_{j+np}\ol f_{j+sp}\int_0^{2\pi}e^{-i(n-s)x}{dx\/2\pi}\rt)=
 \sum_{j=1}^p\sum_{n=-\iy}^{\iy}\|f_{j+np}\|^2=\|f\|^2.
$$
Then $U$ is well defined and has a unique extension to an isometry.
To see that $U$ is onto $\mH$ we compute $U^*$.
We define
$$
 (U^*g)_{j+pn}=\int_0^{2\pi}e^{inx}g_j(x){dx\/2\pi}\in\C^m,\ \
 (j,n)\in\N_{p}\ts\in\Z.
$$
A direct computation shows that this is indeed the formula for
the adjoint of $U$. Moreover,
$$
\|U^*g\|^2={1\/(2\pi)^2}\sum_{j=1}^p\sum_{n=-\iy}^{\iy}\|(U^*g)_{j+pn}\|^2=\sum_{j=1}^p\sum_{n=-\iy}^\iy\lt\|\int_0^{2\pi}e^{inx}g_j(x){dx\/2\pi}\rt\|^2
$$
$$
 =\sum_{j=1}^p\int_0^{2\pi}\|g_j(x)\|^2{dx\/2\pi}=
 \int_0^{2\pi}\sum_{j=1}^p\|g_x(j)\|^2{dx\/2\pi}=\|g\|^2,
$$
where we have used the Parseval relation for the Fourier
series.

We verify \er{5004}. Using $a_{n+p}=a_n,b_{n+p}=b_n$
 we obtain for $f\in \mS(\Z)$
$$
(U\cJ f)_j(x)=\sum_{n=-\iy}^{\iy}e^{-inx}(Jf)_{j+pn}
=\sum_{n=-\iy}^{\iy}e^{-inx} (a_{j+pn-1}f_{j+pn-1}+b_{j+pn}f_{j+pn}+a_{j+pn}f_{j+pn+1})
$$
$$ %=\sum_{n=-\iy}^{\iy}e^{-inx}(a_{j-1}f_{j+pn-1}+b_{j}f_{j+pn}+a_{j}f_{j+pn%+1})=
=a_{j-1}\sum_{n=-\iy}^{\iy}e^{-inx}f_{j+pn-1}+b_j
 \sum_{n=-\iy}^{\iy}e^{-inx}f_{j+pn}+
+a_j\sum_{n=-\infty}^{\infty}e^{-inx}f_{j+pn+1}
$$$$=a_{j-1}(Uf)_x(j-1)+b_j(Uf)_x(j)+a_j(Uf)_x(j),\
 \ \ j\in\N_p,
$$
where we define $(Uf)_x(0)$ и $(Uf)_x(p+1)$ (since
$(Uf))x(j)$ is defined for $j\in\N_p$) by
$$
 (Uf)_x(0)=\sum_{n=-\infty}^{\infty}e^{-inx}f_{pn}=
 e^{-ix}\sum_{n=-\infty}^{\infty}e^{-i(n-1)x}f_{p+p(n-1)}=e^{-ix}(Uf)_x(p),
$$
$$
 (Uf)_x(p+1)=\sum_{n=-\infty}^{\infty}e^{-inx}f_{p+1+pn}=
 e^{ix}\sum_{n=-\infty}^{\infty}e^{-i(n+1)x}f_{1+p(n+1)}=e^{ix}(Uf)_x(1),
$$
which yields $(U\cJ f)(x)=L(e^{ix})(Uf)(x)$.

The eigenvalues $\l_n(e^{ix}), n\in \N_{mp}$ are piecewise real analytic functions on $[0,2\pi]$ , see [RS]. This yields $\s_{sc}(\cJ)=\es$.
Moreover,  the standard simple arguments from \cite{DS} yield  that the operator  $\cJ$ has hot have eigenvalues. Then standard arguments
from \cite{RS} yield \er{Tdi-3}.
\BBox

\begin{lemma} The following identity hold true
\[
 \lb{3001}
 \det(L(\t)-zI_{pm})={D(z,\t)\/c\t^{m}}
 =(-1)^{pm}\prod_{n=1}^{pm}(z-\l_n(\t)),\ \ all \qq \t\not=0,z.
\]
\end{lemma}
\no{\bf Proof.}  Let $L(\t)f=\l(\t)f$ for some eigenvalue
$\l(\t)$ and  some eigenfunction $f=(f_n)_{n=1}^p\in \C^{pm}$.
If $f_0=\t^{-1}f_p$, then the definition of the matrix  $\mM_p$ gives
$ \mM_p(\l(\t))(f_0,f_1)^\top=\t(f_0, f_1)^\top$.
Thus $\t$ is a multiplier of $\mM_p(\l(\t))$ and
$\l(\t)$ is a zero of $D(\cdot,\t)$.

Firstly, let all eigenvalues $\l_n(\t), n\in \N_{mp}$ of $L(\t)$ be distinct for some $\t\in \C$.
Then $\l_n(\t), n\in \N_{mp}$ are zeros of $D(\cdot,\t)$,
which yields \er{3001} for all $\t,z\in\C$, since the orders of $\det(L(\t)-z I_{mp})$ and $D(z,\t)$  are coincide.

Secondly, consider the general case. Let $r_n=n, n\in \N_{mp}$ and let
$r_{n+pm}=r_{n}$ for all $n\in \Z$.  Define the operator Jacobi
$\cJ_t=\cJ+t\diag(r_n)_{n\in\Z}, t\in\R$
 and the corresponding matrix
$L_t(\t)=L(\t)+t\diag(r_n)_{n=1}^{pm}$ (given by \er{5004})
and let $D_t(z,\t)$ be the corresponding determinant.
Then all eigenvalues $\l_{t,n}(\t), n\in \N_{mp}$ of $L_t(\t)$ are distinct for some $\t\in \C$ and $t\to\iy$, since $\l_{t,n}(\t)=tr_n+o(t)$ as $t\to\iy$.  Thus we obtain
\[
\lb{f1}
 \det(L_t(\t)-z I_{pm})={D_t(z,\t)\/c\t^{m}}=(-1)^{pm}\prod_{n=1}^{pm}(z-\l_{t,n}(\t)),
\]
for some $\t\ne 0$ and for all $z\in\C$ and all large $t$.
The functions in \er{f1} are polynomials in $t$, then
identities \er{f1} hold true for all $t\in \R$
and in particular, at $t=1$,
which yields \er{3001}. \BBox

 \no {\bf Acknowledgments.}
Evgeny Korotyaev was partly supported by DFG project BR691/23-1.
The some part of this paper was written at the Math. Institute of
Humboldt Univ., Berlin; Anton Kutsenko is grateful to the
Institute for the hospitality.

\end{document}